\documentclass[english]{elsarticle}
\usepackage[T1]{fontenc}
\usepackage[latin9]{inputenc}
\usepackage{color}
\usepackage{float}
\usepackage{amsmath}
\usepackage{amssymb}
\usepackage{graphicx}

\makeatletter
\usepackage{fullpage}
\usepackage{lmodern}
\usepackage{url}
\usepackage{mathtools,amssymb}

\makeatother

\usepackage{babel}
\begin{document}
\title{A mass, momentum, and energy conservative dynamical low-rank scheme for the Vlasov equation}
\author[uibk]{Lukas Einkemmer\corref{cor1}} \ead{lukas.einkemmer@uibk.ac.at}
\author[llnl]{Ilon Joseph}
\address[uibk]{Department of Mathematics, University of Innsbruck, Austria}
\address[llnl]{Physics Division, Lawrence Livermore National Laboratory, California, USA}
\cortext[cor1]{Corresponding author}
\begin{abstract} The primary challenge in solving kinetic equations, such as the Vlasov equation, is the high-dimensional phase space. In this context, dynamical low-rank approximations have emerged as a promising way to reduce the high computational cost imposed by such problems. However, a major disadvantage of this approach is that the physical structure of the underlying problem is not preserved. In this paper, we propose a dynamical low-rank algorithm that conserves mass, momentum, and energy as well as the corresponding continuity equations. We also show how this approach can be combined with a conservative time and space discretization. 
\end{abstract}  
\maketitle

\section{Introduction}

Solving kinetic equations efficiently is important in applications
ranging from plasma physics to radiative transfer. The main challenge
in this context is the up to six-dimensional phase space and the associated
unfavorable scaling of computational cost and memory requirements,
usually referred to as the curse of dimensionality. Because of this,
particle methods, such as the particle in cell (PIC) scheme, have
been and are still widely used. However, it also well known that particle
methods miss or do not resolve certain physical phenomena (see, e.g.,
\citep{camporeale2016}) or require an immense number of particles
thereby negating their advantage. In addition, complexity reduction
techniques such as sparse grids have been investigated. While they
can provide some advantage compared to a full grid simulation the
gain is usually modest \citep{Kormann2014}.

More recently, using dynamical low-rank approximations to solve kinetic
problems has received considerable interest. Such methods have been
developed for both the Vlasov equation \citep{Kormann15,El18,EL18_cons,einkemmer2020low}
and radiation transport problems \citep{ding2019error,peng2019low,einkemmer2020asymptotic}.
Dynamical low-rank integrators approximate a six-dimensional Vlasov
equation by a set of only three dimensional advection problems. Moreover,
they have a range of properties that makes them well suited for performing
kinetic simulations. In particular, such methods are able to resolve
filamentation and are almost exact if the dynamics can be well represented
by a linearized equation \citep{einkemmer2020low}. In addition, dynamical
low-rank schemes can capture the limiting fluid or diffusive regime
\citep{E18,ding2019error,einkemmer2020asymptotic,peng2020high}. Therefore,
dynamical low-rank integrators can drastically decrease the numerical
effort that is required to solve a number of kinetic problems. This
enables the simulation of such problems on desktop computers or small
clusters that are otherwise either unfeasible or require large supercomputers.

A major disadvantage of dynamical low-rank integrators, however, is
that they do not respect the physical structure of the underlying
equations. In particular, mass, momentum, and energy are not preserved
by the low-rank approximation. This is in stark contrast to Eulerian
and semi-Lagrangian Vlasov solvers, where at least mass and momentum
are conserved \citep{palmroth2018review} and methods with good long-time
behavior with respect to energy have been obtained \citep{CEF15},
and particle methods, where usually mass and either momentum or energy
can be conserved \citep{verboncoeur2005particle}. 

Some approaches to improve this deficiency have been proposed. In
\citep{peng2019low} the solution is simply rescaled such that the
mass is preserved. Such an approach, however, can not be extended
to simultaneously conserve momentum or energy. It also does not respect
the underlying continuity equation for the mass density, which could
actually be more important for the long time behavior of the integrator
\citep{EL18_cons}. In \citep{EL18_cons} a method based on Lagrange
multipliers is proposed. This approach succeeds in improving the conservative
properties of the low-rank approach, but it does not allow us to simultaneously
conserve both the continuity equations for density and momentum as
well as the corresponding invariants. It also does not change the
dynamical low-rank approximation. Instead, this method adds a correction
to each step of the low-rank integrator. \textcolor{black}{In \citep{peng2020high}
the fluid moments are integrated explicitly. These moments are then
coupled to a low-rank approximation that resolves the kinetic dynamics.
To enforce conservation, a correction is added to the low-rank part
of the algorithm that requires, similar to \citep{EL18_cons}, the
solution of a linear system of equations.}

Neither of the conservative or quasi-conservative methods developed
in the literature solve the fundamental problem. Namely, that the
classic dynamical low-rank approximation does not take the structure
of the equations that we are trying to solve into account. In this
paper we introduce such an approach that conserves both mass, momentum,
and energy and ensures that the corresponding continuity equations
are satisfied as well. The proposed method is based on the observation
that if certain functions of velocity belong to the approximation
space, then the desired conservation follows. To accomplish this two
steps are necessary. First, we have to formulate the dynamical low-rank
scheme in a different function space than the $L^{2}$ space that
is usually used (we will use an appropriately weighted $L^{2}$ space
instead). Second, we have to constrain the low-rank factors such that
the desired functions belong to the approximation space at all times.
This is done by fixing certain basis functions and using a modified
Petrov--Galerkin condition in such a way that this is compatible
with the remainder of the dynamical low-rank approximation. This results
in equations of motions that are somewhat different from the ones
used in the classic dynamical low-rank algorithm introduced in \citep{Koch2010}
(and that subsequently was used heavily in the literature; see, e.g,
\citep{Nonnenmacher2008,musharbash2015error,El18,EL18_cons,peng2019low,einkemmer2020low,einkemmer2020asymptotic}).
We also introduce a time and space discretization of the resulting
equations of motions that conserves mass and momentum as well as a
discrete version of the corresponding continuity equations up to machine
precision. 

The remainder of the paper is structured as follows. In section \ref{sec:vp-conservation}
we introduce the Vlasov--Poisson equation and discuss the invariants
and corresponding physical structure of this model. This is followed
by a recollection of the classic dynamical low-rank integrator in
section \ref{sec:classic-dlr}. The proposed conservative dynamical
low-rank integrator is then introduced in section \ref{sec:cons-dlr}
and its properties are discussed in section \ref{sec:Conservation}.
The fully discretized integrator, i.e.~time and space discretization,
is discussed in section \ref{subsec:Time-integration}. Finally, a
number of numerical simulations are presented in section \ref{sec:numerical-experiments}.

\section{Vlasov--Poisson equations and conservation\textcolor{black}{\label{sec:vp-conservation}}}

In this work we consider the Vlasov--Poisson equations
\begin{align}
 & \partial_{t}f(t,x,v)+v\cdot\nabla_{x}f(t,x,v)-E(f)(t,x)\cdot\nabla_{v}f(t,x,v)=0,\label{eq:vlasov}\\
 & \nabla_{x}\cdot E(f)(t,x)=1-\int f(t,x,v)\,\mathrm{d}v,\qquad\;\;\nabla_{x}\times E(f)(t,x)=0\label{eq:poisson}
\end{align}

on $(x,v)\in\Omega_{x}\times\Omega_{v}$ with $\Omega_{x}\subset\mathbb{R}^{d}$
(physical space) and $\Omega_{v}\subset\mathbb{R}^{d}$ (velocity
space), $d\le3$. The sought after quantities are the particle-density
function $f$ and the electric field $E$. This equation models the
dynamics of electrons that are subject to a homogeneous ion background
distribution. We consider the single species Vlasov--Poisson equation
here, but it should be emphasized that the described algorithm can
be easily extended to the multi-particle case as well as to more complicated
models, such as the Vlasov--Maxwell equations.

The dynamics of the Vlasov--Poisson equations conserves a number
of physically relevant invariants. In particular, mass $M$, momentum
$J$
\[
M(t)=\int_{\Omega}f(t,x,v)\,d(x,v),\qquad J(t)=\int_{\Omega}vf(t,x,v)\,d(x,v)
\]

and energy $\mathcal{E}$ 
\[
\mathcal{E}(t)=\frac{1}{2}\int_{\Omega}v^{2}f(t,x,v)\,d(x,v)+\frac{1}{2}\int_{\Omega_{x}}E(t,x)^{2}\,dx
\]
are invariants of the solution to the continuous equations. We note
that in the single species setting mass conservation is equivalent
to charge conservation (since charge and mass are proportional) and
momentum conservation is equivalent to current conservation.

For each of these invariants, a continuity or moment equation, that
is posed in physical space only, is satisfied by an associated density.
In the case of mass, the mass density
\[
\rho(t,x)=\int_{\Omega_{v}}f(t,x,v)\,dv
\]

satisfies
\[
\partial_{t}\rho(t,x)+\nabla_{x}\cdot j(t,x)=0,\qquad j(t,x)=\int_{\Omega_{v}}vf(t,x,v)\,dv.
\]

Conservation of mass can be easily derived from the continuity equation
by integrating in $x$. It is in fact these continuity equations that
our dynamical low-rank integrator satisfies and from which conservation
of the global invariants follows. We note that this is a much stronger
result than simply conserving the invariants, as the present approach
also preserves the underlying physical structure of the equation.

The associated density for momentum conservation is the momentum density
$j(t,x)$ which satisfies the following continuity equation
\[
\partial_{t}j(t,x)+\nabla_{x}\cdot\sigma(t,x)=-E(t,x)\rho(t,x),\qquad\sigma(t,x)=\int_{\Omega_{v}}(v\otimes v)f(t,x,v)\,dv.
\]

We can derive conservation of momentum by recognizing that $E(1-\rho)=\nabla\cdot(E\otimes E-\tfrac{1}{2}E^{2})$
and integrating in physical space. Note that due to the normalization
of the particle-density function we have $\int E\,dx=0$. The energy
density
\[
e(t,x)=\frac{1}{2}\int_{\Omega_{v}}v^{2}f(t,x,v)\,dv+\frac{1}{2}E^{2}(t,x)
\]
satisfies the following continuity equation
\[
\partial_{t}e(t,x)+\nabla_{x}\cdot Q(t,x)=E(t,x)\cdot(\partial_{t}E(t,x)-j(t,x)),\qquad Q(t,x)=\frac{1}{2}\int_{\Omega_{v}}vv^{2}f(t,x,v)\,dv.
\]

We derive global conservation, i.e.~conservation of energy, by
\[
\partial_{t}\mathcal{E}(t,x)=\partial_{t}\int e(t,x)\,dx=\int E(t,x)\cdot(\partial_{t}E(t,x)-j(t,x))\,dv=0.
\]
The last equality follows since $\partial_{t}E(t,x)=j(t,x)$ is just
the electrostatic version of Ampere's law. A more detailed discussion
can be found in \citep{sonnendruecker_lecturenotes}. Looking at the
conserved quantities in light of their corresponding continuity equations
also relates the kinetic equation considered here to their corresponding
fluid model. For more details we refer the reader to \citep{Hu2017,einkemmer2020asymptotic,E18}.

\section{The classical dynamical low-rank scheme\textcolor{black}{\label{sec:classic-dlr}}}

For a dynamical low-rank scheme we seek an (approximate) solution
of the form 
\begin{equation}
f(t,x,v)=\sum_{i,j=1}^{r}X_{i}(t,x)S_{ij}(t)V_{j}(t,v),\label{eq:lowrank}
\end{equation}
where $r$ is the rank of the approximation. The dynamics is represented
by the low-rank factors $X_{i}$ and $V_{j}$, that only depend on
physical space $x$ and velocity space $v$, respectively, and the
low-rank factor $S_{ij}$, which carries no spatial or velocity dependence.

From a physical point of view, the decomposition into physical space
and velocity space is natural. If the dynamics is integrable then
there are conserved action variables $j(x,v)$ and conjugate angle
variables $\theta(x,v)$. In this case, the natural basis functions
are $\Theta_{k}(\theta,t)=e^{ik\theta}$ and $J_{k}(j,t)=\delta(j-j_{k})$.
Note, however, that there is a nonlinear coordinate transformation
that must be determined in order for this representation to hold.
The separation into physical space and velocity used here does hold
for linearized equations, such as the Case-Van Kampen eigenfunctions
\citep{van1955theory} for the linear Vlasov--Poisson equation. For
a mathematical perspective in the context of the dynamical low-rank
approximation, see \citep{einkemmer2020low}. In the collisional case,
the nonlinear collision operator generates scattering in $v$ and
is considered local in $x$. This simplification would not hold in
generalized phase space coordinates $\{\theta,j\}$. Hence, when collisions
are dominant, the decomposition between physical space and velocity
space is still the natural one. Due to the fact that sharp structures,
such as shocks and boundary layers, can form in the $x$ direction,
good basis functions are localized in $x$. For the $v$ direction,
good basis functions include functions orthogonal over a Maxwellian
distribution and eigenfunctions of the collision operator. The low-rank
approximation can accommodate this without any difficulty, as has
been demonstrated in \citep{E18,ding2019error,einkemmer2020asymptotic}.

From now on, we will, for reasons of simplicity, write $f$ for $f(t,x,v)$,
$X_{i}$ for $X_{i}(t,x)$, etc. The equations of motions for the
low-rank factors $X_{i}$, $S_{ij}$, and $V_{j}$ are then determined
by imposing a Galerkin condition, see \citep{Koch2010}. Following
this argument for the PDE case we obtain \citep{El18} 
\begin{align}
\partial_{t}\left(\sum_{i}X_{i}S_{ij}\right) & =(V_{j}D[f])_{v},\label{eq:classic-lr-K}\\
\partial_{t}S_{ij} & =(X_{i}V_{j}D[f])_{xv},\label{eq:classic-lr-S}\\
\partial_{t}\left(\sum_{j}S_{ij}V_{j}\right) & =(X_{i}D[f])_{x},\label{eq:classic-lr-L}
\end{align}
where
\[
(f)_{x}=\int_{\Omega_{x}}f\,dx,\qquad(f)_{v}=\int_{\Omega_{v}}f\,dv,\qquad(f)_{xv}=\int_{\Omega}f\,d(x,v).
\]
For the Vlasov--Poisson equation, $D[f]=-v\cdot\nabla_{x}f+E\cdot\nabla_{v}f$
with $f$ in the form of equation (\ref{eq:lowrank}). We note that
the low-rank factors $X_{i}$ and $V_{j}$ are orthonormal; that is,
they satisfy $(X_{i}X_{j})_{x}=\delta_{ij}$ and $(V_{i}V_{j})_{v}=\delta_{ij}$.
The main utility of the low-rank approximation is that instead of
a problem in dimension $2d$, as in equation (\ref{eq:vlasov}), we
only have to solve $2r$ equations of dimension $d$ and $r$ ordinary
differential equations. This is the reason why the dynamical low-rank
approximation drastically reduces the memory and computational effort
required to solve such problems, assuming that $r$ is not too large,
which is often true in practice. For a more detailed discussion we
refer the reader to \citep{El18,einkemmer2020low} and the subsequent
discussion in section \ref{sec:cons-dlr}.

Equations (\ref{eq:classic-lr-K})-(\ref{eq:classic-lr-L}), however,
do not preserve the invariants of the original system. That is, even
if the system is solved exactly (i.e.~no time or space error is introduced)
conservation of mass, momentum, and energy is lost. This is not very
surprising as the dynamical low-rank approximation described above
does not take the physical structure of the equations into account.
Thus, there is no mechanism that prevents, e.g., mass to be removed
from the system due to the truncation performed by the low-rank approximation. 

\section{The conservative dynamical low-rank scheme\label{sec:cons-dlr}}

In this section, we propose a novel dynamical low-rank integrator
with equations of motions, that will take the place of (\ref{eq:classic-lr-K})-(\ref{eq:classic-lr-L}),
that ensure mass, momentum, and energy conservation.

For the classical dynamical low-rank integrator we have
\[
\partial_{t}f=\sum_{ij}\partial_{t}(X_{i}S_{ij})V_{j}-\sum_{ij}X_{i}(\partial_{t}S_{ij})V_{j}+\sum_{ij}X_{i}\partial_{t}(S_{ij}V_{j})
\]
and thus by integrating over $v$
\begin{equation}
\partial_{t}\rho=\sum_{j}(V_{j}D[f])_{v}(1V_{j})_{v}-\sum_{ij}X_{i}(X_{i}V_{j}D[f])_{xv}(1V_{j})_{v}+\sum_{i}X_{i}(X_{i}D[f])_{xv},\label{eq:rhoconsd}
\end{equation}

where we have used equations (\ref{eq:classic-lr-K})-(\ref{eq:classic-lr-L}).
Now, if we could ensure that $1\in\overline{V}$, where $\overline{V}=\text{span}\{V_{j}\}$,
then we could find coefficients $c_{j}$ such that $1=\sum_{j}c_{j}V_{j}$.
Plugging this into equation (\ref{eq:rhoconsd}) we get
\[
\partial_{t}\rho=\Bigl(D[f]\sum_{j}c_{j}V_{j}\Bigr)_{v}-\sum_{i}X_{i}\Bigl(X_{i}D[f]\sum_{j}c_{j}V_{j}\Bigr)_{xv}+\sum_{i}X_{i}(X_{i}D[f])_{xv}
\]
and thus
\[
\partial_{t}\rho=\left(D[f]\right)_{v}-\sum_{i}X_{i}(X_{i}D[f])_{xv}+\sum_{i}X_{i}(X_{i}D[f])_{xv}=(D[f])_{v}.
\]

By using $D[f]=-v\cdot\nabla_{x}f+E\cdot\nabla_{v}f$ we at once obtain
the continuity equation
\[
\partial_{t}\rho+\nabla\cdot j=0,
\]

which would imply mass conservation. However, $1\in\overline{V}$
is clearly not true as constant functions do not lie in $L^{2}(\mathbb{R}^{d})$.
Thus, the continuity equation is not satisfied and the classic dynamical
low-rank integrator does not conserve mass. A very similar argument
can be made for the momentum, with $v\in\overline{V}$, and energy,
with $v^{2}\in\overline{V}$. One might object at this point that
in a numerical simulation we necessarily use a truncated domain and
thus the functions $1$, $v,$ and $v^{2}$ lie in $L^{2}(\Omega_{v})$.
However, on a finite domain the approximation spaces have to be equipped
with appropriate boundary conditions (usually either periodic or homogeneous
Dirichlet conditions are imposed). \textcolor{black}{The functions
$1$, $v$, and $v^{2}$ do not satisfy boundary conditions that are
compatible with the problem. In the periodic case it is clear that
$v$ and $v^{2}$ are not periodic and assuming so would incur a large
numerical error. In the case of homogenous Dirichlet boundary conditions
$1$, $v$, $v^{2}$ are not zero at the boundary and thus do not
lie in the desired approximation space.}

To obtain a conservative dynamical low-rank integrator we proceed
as follows. First, we will use a weighted function space that includes
the possibility of representing constant functions in $\overline{V}.$
However, this on its own is not yet sufficient. In fact, the dynamic
low-rank algorithm automatically chooses appropriate basis functions
in order to minimize the overall error, according to some Galerkin
condition, in the particle-density function $f$. There is no guarantee
that such a choice satisfies $1\in\text{span}\{\overline{V}\}.$ Thus,
we have to constrain the approximation space in such a way that constant
functions are always part of the basis. This, in particular, requires
us to modify the Galerkin condition that is used in the classic integrator.
Those two ideas combined allow us to formulate a conservative dynamical
low-rank algorithm. The details of this procedure will be the topic
of the remainder of this section.

For the conservative dynamical low-rank scheme the solution is approximated
by a function of the following form
\begin{equation}
f(t,x,v)=f_{0}(t,x,v)\sum_{i,j=1}^{r}X_{i}(t,x)S_{ij}(t)V_{j}(t,v),\label{eq:conslr}
\end{equation}
where $r$ is the rank of the approximation and $f_{0}(x,v)=f_{0x}(t,x)f_{0v}(t,v)$
is a, yet to be chosen, weight function. The low-rank factors $X_{i}(t,\cdot)$
and $V_{j}(t,\cdot)$ are assumed to lie in the $L^{2}$ spaces weighted
by $f_{0x}$ and $f_{0v}$, respectively. That is, $X_{i}(t,\cdot)\in L^{2}(\Omega_{x},f_{0x})$
and $V_{j}(t,\cdot)\in L^{2}(\Omega_{v},f_{0v})$. The corresponding
weighted inner products are denoted by
\[
\langle X_{i}(t,\cdot),X_{j}(t,\cdot)\rangle_{x}=\int_{\Omega_{x}}f_{0x}(t,x)X_{i}(t,x)X_{j}(t,x)\,dx
\]

and 
\[
\langle V_{i}(t,\cdot),V_{j}(t,\cdot)\rangle_{v}=\int_{\Omega_{v}}f_{0v}(t,v)V_{i}(t,v)V_{j}(t,v)\,dv,
\]
respectively.

In the following we will choose $f_{0x}(t,x)=1$ and $f_{0v}(t,v)=f_{0v}(v)$,
i.e. $f_{0v}$ is time independent. The choice of $f_{0v}$ must guarantee
that $1$, $v$ , and $v^{2}$ lie in $L^{2}(\Omega_{v},f_{0v})$
in order to obtain conservation of mass, momentum, and energy, respectively.
It is also important to choose the temperature large enough such that
the relevant part of the phase space is captured. For many problems
that model plasma instabilities $f_{0v}(v)=\exp(-v^{2}/2)$ is a good
choice as this represents the equilibrium distribution and guarantees
that all powers of $v$ lie in the approximation space. However, other
than the constraints outlined above the choice of $f_{0v}$ is arbitrary.

The second crucial ingredient of the algorithm is that some of the
functions $V_{j}$ are held fixed as the system evolves in time. We
write
\[
U_{a}(v)=V_{a}(v),\quad1\leq a\leq m\qquad\text{and}\qquad W_{p}(t,v)=V_{p}(t,v),\quad m<p\leq r,
\]
where the $U_{a}(v)$ are fixed, i.e.~they are not changed by the
dynamical low-rank integrator, and the $W_{p}(t,v)$ are allowed to
vary in time according to the low-rank algorithm. In the following,
indices $i$, $j$, $k$, $l$ span $1,\dots,r$, indices $p$, $q$
span $m+1,\dots,r$, and the indices $a$, $b$ span $1,\dots,m$.
For example, to obtain mass, momentum, and energy conservation in
1+1 dimension and for $f_{0v}(v)=\exp(-v^{2}/2)$ we choose $m=3$
with $U_{1}(v)=1/\Vert1\Vert$, $U_{2}(v)=v/\Vert v\Vert$, and $U_{3}(v)=(v^{2}-1)/\Vert v^{2}-1\Vert$.
It is clear that then the orthogonality constraint $\langle U_{a},U_{b}\rangle_{v}=\delta_{ab}$
is satisfied. As usual, we also impose the orthogonality conditions
$\langle W_{p},V_{j}\rangle=\delta_{pj}$. Note that this implies
that the $W_{p}$ are orthogonal both with respect to each other and
with respect to the $U_{a}$.

The equations of motion for the low-rank factors are derived by considering
the low-rank manifold of all functions of a given rank $r$. In our
case this manifold $\mathcal{M}$ can be written as
\begin{align*}
\mathcal{M} & =\biggl\{ f\in L^{2}(\Omega,f_{0})\colon f(x,v)=f_{0v}(v)\sum_{ij}X_{i}(x)S_{ij}V_{j}(v)\text{ with invertible }S=(S_{ij})\in\mathbb{R}^{r\times r},\\
 & \qquad X_{i}\in L^{2}(\Omega_{x}),\,V_{j}\in L^{2}(\Omega_{v},f_{0v})\text{ with }\langle X_{i},X_{j}\rangle_{x}=\delta_{ij},\,\langle V_{i},V_{j}\rangle_{v}=\delta_{ij}\biggr\}.
\end{align*}
As usual we will impose the gauge conditions $\langle\partial_{t}X_{i},X_{j}\rangle_{x}=0$
and $\langle\partial_{t}W_{p},W_{q}\rangle_{v}=0$ in order to make
sure that the dynamics of the low-rank factors is uniquely determined.
The tangent space of the manifold is then
\begin{align*}
\mathcal{T}_{f}\mathcal{M} & =\biggl\{\dot{f}\in L^{2}(\Omega,f_{0})\colon\dot{f}=f_{0v}\sum_{ij}\left(\dot{X}_{i}S_{ij}V_{j}+X_{i}\dot{S}_{ij}V_{j}\right)+f_{0v}\sum_{ip}X_{i}S_{ip}\dot{W}_{p},\\
 & \qquad\text{with }\dot{S}\in\mathbb{R}^{r\times r},\,\dot{X}_{i}\in L^{2}(\Omega_{x}),\,\dot{V}_{j}\in L^{2}(\Omega_{v},f_{0v}),\text{ and }\langle X_{i},\dot{X}_{j}\rangle_{x}=0,\ \langle V_{i},\dot{W_{q}}\rangle_{v}=0\biggr\},
\end{align*}
where we have used the fact that $\partial_{t}U_{a}=0$. While this
construction is similar to the classic dynamical low-rank integrator,
it is crucial that we are cognizant for which low-rank factors the
gauge conditions are imposed and that we make sure that the entire
set of the $V_{j}$ are orthogonal to each other (and not only the
$W_{p}$). 

We now derive the equations of motion for the low-rank factors $X_{i}$,
$S_{ij}$, and $V_{ij}$ that satisfy the following Petrov--Galerkin
condition
\begin{equation}
\left(\frac{{\color{black}\nu}}{f_{0v}},(\partial_{t}f-D[f])\right)_{xv}=0\qquad\forall{\color{black}\nu}\in T_{f}\mathcal{M},\label{eq:galerkin-condition}
\end{equation}

where $(f,g)_{xv}=(fg)_{xv}$ and $\partial_{t}f\in\mathcal{T}_{f}\mathcal{M}$.
\textcolor{black}{Note that $(\cdot,\cdot)_{x}$, $(\cdot,\cdot)_{v}$,
and $(\cdot,\cdot)_{xv}$ denote the usual $L^{2}$ inner products,
while $\langle\cdot,\cdot\rangle_{x}$, $\langle\cdot,\cdot\rangle_{v}$,
and $\langle\cdot,\cdot\rangle_{xv}$ denote the $L^{2}$ inner products
weighted with $f_{0x}$, $f_{0v}$, and $f_{0}$, respectively. Since
we have assumed $f_{0x}(t,x)=1$, it holds that $(\cdot,\cdot)_{x}=\langle\cdot,\cdot\rangle_{x}$.}

The Petrov--Galerkin condition is different from the classic dynamical
low-rank scheme; an additional factor of $1/f_{0v}$ has been introduced.
The reason for this is that in the following our goal is to choose
an appropriate $v\in T_{f}\mathcal{M}$ such that the equation of
motion for a specific low-rank factor is isolated. This change in
the Petrov--Galerkin condition is what makes this possible in the
weighted approximation spaces that we consider here.

\textbf{Equations for $X_{i}$: }In this case there is little difference
to the classic algorithm. We consider a family of test functions $\nu_{k}=f_{0v}\chi(x)V_{k}$,
where $\chi$ is an arbitrary function of $x$. Since we can write
$\nu_{k}=f_{0v}\sum_{ij}\dot{X_{i}}S_{ij}V_{j}$ with $\dot{X}_{i}=\chi(x)S_{ki}^{-1}$
it holds that $\nu_{k}\in T_{f}\mathcal{M}$. The Petrov--Galerkin
condition (\ref{eq:galerkin-condition}) then becomes
\[
\biggl(V_{k}\chi(x),f_{0v}\sum_{ij}\left(\dot{X}_{i}S_{ij}V_{j}+X_{i}\dot{S}_{ij}V_{j}\right)+f_{0v}\sum_{ip}X_{i}S_{ip}\dot{W}_{p}\biggr)_{xv}=\left(V_{k}\chi(x),D[f]\right)_{xv}.
\]
We can rewrite this as
\[
\biggl\langle V_{k}\chi(x),\sum_{ij}\left(\dot{X}_{i}S_{ij}V_{j}+X_{i}\dot{S}_{ij}V_{j}\right)+\sum_{ip}X_{i}S_{ip}\dot{W}_{p}\biggr\rangle_{xv}=\left(V_{k}\chi(x),D[f]\right)_{xv}
\]

and thus by using the orthogonality and gauge conditions as well as
the fact that $\chi$ is arbitrary we obtain
\begin{equation}
\sum_{i}\dot{X}_{i}S_{ik}=\left(V_{k},D[f]\right)_{v}-\sum_{i}X_{i}\dot{S}_{ik},\label{eq:evol-X}
\end{equation}
which is precisely the first equation of motion that is obtained for
the classic algorithm (see, e.g., \citep{El18}). 

Now, we can plug the right hand-side of the Vlasov--Poisson equation
into the right-hand side of equation (\ref{eq:evol-X}). This yields
\begin{align*}
\left(V_{k},D[f]\right)_{v} & =\left(V_{k},-v\cdot\nabla_{x}f+E\cdot\nabla_{v}f\right)_{v}\\
 & =-\sum_{ij}\left(V_{k},f_{0v}v\cdot(\nabla_{x}X_{i})S_{ij}V_{j}\right)_{v}+\sum_{ij}\left(V_{k},X_{i}S_{ij}E\cdot\nabla_{v}(f_{0v}V_{j})\right)_{v}\\
 & =-\sum_{ij}\langle V_{k},vV_{j}\rangle_{v}\nabla_{x}X_{i}S_{ij}+\sum_{ij}E\cdot(V_{k},\nabla_{v}(f_{0v}V_{j}))_{v}X_{i}S_{ij}\\
 & =-\sum_{ij}c_{kj}^{1}\nabla_{x}X_{i}S_{ij}+\sum_{ij}(c_{kj}^{2}\cdot E)X_{i}S_{ij},
\end{align*}
with
\[
c_{kj}^{1}=\langle V_{k},vV_{j}\rangle_{v},\qquad\qquad c_{kj}^{2}=(V_{k},\nabla_{v}(f_{0v}V_{j}))_{v}.
\]

We note that due to the weight function $f_{0v}$ the coefficients
$c_{kl}^{1}$ and $c_{kl}^{2}$ are changed compared to the classic
algorithm in \citep{El18}.

\textbf{Equations for $W_{p}$: }In this case we consider the family
of test functions $\nu_{q}=f_{0v}\zeta(v)\sum_{i}X_{i}S_{iq}$, where
$\zeta$ is an arbitrary function of $v$. Those $\nu_{q}$ lie in
the tangent space as we can write $\nu_{q}=f_{0v}\sum_{ip}X_{i}S_{ip}\dot{W}_{p}$
with $\dot{W}_{p}=\delta_{pq}\zeta(v)$. The Petrov--Galerkin condition
(\ref{eq:galerkin-condition}) then becomes
\[
\sum_{i}\biggl(\zeta(v)X_{i}S_{iq},f_{0v}\sum_{kl}\left(\dot{X}_{k}S_{kl}V_{l}+X_{k}\dot{S}_{kl}V_{l}\right)+f_{0v}\sum_{kp}X_{k}S_{kp}\dot{W}_{p}\biggr)_{xv}=\sum_{i}(\zeta(v)X_{i}S_{iq},D[f])_{xv}
\]

and thus (since $\zeta$ is arbitrary and using the gauge conditions)
\begin{equation}
\sum_{ip}S_{iq}S_{ip}\dot{W}_{p}=\frac{1}{f_{0v}}\sum_{i}S_{iq}(X_{i},D[f])_{x}-\sum_{il}S_{iq}\dot{S}_{il}V_{l}.\label{eq:evol-W}
\end{equation}

On the left-hand side we have the matrix $T_{qp}=\sum_{i}S_{iq}S_{ip}$.
Since $S$ has full rank the same is true for $T$ and thus we can
invert it in order to obtain the equations of motion for the $W_{p}$.

For the Vlasov--Poisson equation we have
\begin{align*}
\frac{1}{f_{0v}}(X_{i},D[f])_{x} & =\frac{1}{f_{0v}}\left(X_{i},-v\cdot\nabla_{x}f+E\cdot\nabla_{v}f\right)_{x}\\
 & =-\frac{1}{f_{0v}}\sum_{kl}f_{0v}\langle X_{i},\nabla_{x}X_{k}\rangle_{x}\cdot vS_{kl}V_{l}+\frac{1}{f_{0v}}\sum_{kl}S_{kl}\nabla_{v}(f_{0v}V_{l})\cdot\langle X_{i},EX_{k}\rangle_{x}\\
 & =-\sum_{kl}(v\cdot d_{ik}^{2})S_{kl}V_{l}+\frac{1}{f_{0v}}\sum_{kl}d_{ik}^{1}[E]\cdot\nabla_{v}(f_{0v}S_{kl}V_{l})\\
 & =-\sum_{kl}(v\cdot d_{ik}^{2})S_{kl}V_{l}+\sum_{kl}d_{ik}^{1}[E]\cdot\left[\nabla_{v}(S_{kl}V_{l})+\nabla_{v}(\log f_{0v})S_{kl}V_{l}\right],
\end{align*}
where 
\[
d_{ik}^{1}[E]=\langle X_{i},EX_{k}\rangle_{x},\qquad\qquad d_{ik}^{2}=\langle X_{i},\nabla_{x}X_{k}\rangle_{x}.
\]

\textbf{Equation for $S_{ij}$: }In this case we consider $\nu_{kl}=f_{0v}X_{k}V_{l}$,
which lies in the tangent space as $\nu_{kl}=f_{0v}\sum_{ij}X_{i}\dot{S}_{ij}V_{j}$
with $\dot{S}_{ij}=\delta_{ki}\delta_{jl}$. The Petrov--Galerkin
condition (\ref{eq:galerkin-condition}) becomes
\[
\biggl(X_{k}V_{l},f_{0v}\sum_{ij}\left(\dot{X}_{i}S_{ij}V_{j}+X_{i}\dot{S}_{ij}V_{j}\right)+f_{0v}\sum_{ip}X_{i}S_{ip}\dot{W}_{p}\biggr)_{xv}=\left(X_{k}V_{l},D[f]\right)_{xv}
\]
and thus
\begin{equation}
\dot{S}_{kl}=\left(X_{k}V_{l},D[f]\right)_{xv}.\label{eq:evol-S}
\end{equation}
For the Vlasov--Poisson equation 
\begin{align*}
\left(X_{k}V_{l},D[f]\right)_{xv} & =(X_{k}V_{l},-v\cdot\nabla_{x}f+E\cdot\nabla_{v}f)_{xv}\\
 & =-\sum_{ij}\langle X_{k},\nabla_{x}X_{i}\rangle_{x}\cdot\langle V_{l},vV_{j}\rangle_{v}S_{ij}+\sum_{ij}\langle X_{k},EX_{i}\rangle_{x}\cdot(V_{l},\nabla_{v}(f_{0v}V_{j}))_{v}S_{ij}\\
 & =-\sum_{ij}(d_{ki}^{2}\cdot c_{lj}^{1})S_{ij}+\sum_{ij}(d_{ki}^{1}\cdot c_{lj}^{2})S_{ij}.
\end{align*}
Together equations (\ref{eq:evol-X})-(\ref{eq:evol-S}) are the equations
of motions for the proposed conservative dynamical low-rank integrator.
They take the place of equations (\ref{eq:classic-lr-K})-(\ref{eq:classic-lr-L})
that have been used in the classic algorithm. The primary difference
lies in the equation for the $W_{p}$, where it has to be ensured
that the $W_{p}$ are updated in such a way that they remain orthogonal
to not only the other $W_{p}$ but also to the fixed $U_{a}$ \textcolor{black}{(this
and the desire to keep the $U_{a}$ fixed is the reason why we can
not simply choose $\nu_{q}=\chi(v)X_{q}$ in deriving the equations
for the $W_{p}$ above)}. In addition, due to the use of the weighted
$L^{2}$ spaces all coefficients are changed as well and $f_{0v}$
and its derivatives make an appearance in the equations of motion\@.

\section{Mass, momentum, and energy conservation for the proposed low-rank
approximation\label{sec:Conservation}}

In this section we will show that the numerical scheme derived in
section \ref{sec:cons-dlr} is indeed mass, momentum, and energy conservative. 

To ensure mass conservation we choose, as is discussed in the previous
section, $U_{1}\propto1$. Then, using $K_{j}=\sum_{i}X_{i}S_{ij}$,
the density is given by
\[
\rho=\frac{1}{U_{1}}K_{1}.
\]

This can be easily seen as
\[
\rho=\int_{\Omega_{v}}f\,dv=\frac{1}{U_{1}}\sum_{ij}X_{i}S_{ij}\langle U_{1},V_{j}\rangle_{v}.
\]

Using the orthogonality condition $\langle U_{1},V_{j}\rangle=\delta_{1j}$
we obtain the desired relation.

Deriving the corresponding continuity equation is now straightforward.
We have
\[
\partial_{t}\rho=\frac{1}{U_{1}}\partial_{t}K_{1}.
\]

From equation (\ref{eq:evol-X}) we derive $\partial_{t}K_{1}=\left(V_{1},D[f]\right)_{v}$
and thus
\[
\partial_{t}\rho=\frac{1}{U_{1}}\left(U_{1},D[f]\right)_{v}=\int_{\Omega_{v}}D[f]\,dv.
\]
This is precisely the relation that is obtained for the continuous
evolution, i.e.~the solution of the Vlasov--Poisson equation without
any low-rank approximation present. Since
\begin{align*}
\int D[f]\,dv & =-\nabla_{x}\cdot\int_{\Omega_{v}}vf\,dv+E\cdot\int_{\Omega_{v}}\nabla_{v}f\,dv\\
 & =-\nabla_{x}\cdot j,
\end{align*}
we have
\[
\partial_{t}\rho+\nabla_{x}\cdot j=0.
\]
From the continuity equation conservation of mass follows at once
(by integrating in space).  

For conservation of momentum we proceed in precisely the same manner.
For simplicity we only consider the 1+1 dimensional case here. By
choosing $U_{2}$ such that $v=\Vert v\Vert U_{2}$, i.e. $U_{2}\propto v$,
momentum is given by\textbf{
\[
j=\Vert v\Vert K_{2}.
\]
}

We thus have
\[
\partial_{t}j=\Vert v\Vert\partial_{t}K_{2}=\int_{\Omega_{v}}vD[f]\,dv=-\nabla_{x}\cdot\sigma-E\rho,
\]
which is the desired continuity equation that implies conservation
of momentum. The extension to the multi-dimensional case is straightforward.
We have to use one fixed basis function for each of the directions
in which we want to conserve momentum (e.g.~$U_{2}\propto v_{1}$,
$U_{3}\propto v_{2}$, and $U_{4}\propto v_{3}$).

For energy conservation we can not choose $U_{3}\propto v^{2}$ since
$v^{2}$ is not orthogonal to $1$. Thus, we use $U_{3}\propto v^{2}-1.$
In this setting we have $\langle U_{3},U_{2}\rangle_{v}=\langle U_{3},U_{1}\rangle_{v}=0$,
as desired. We then can represent $v^{2}$ as follows
\[
v^{2}=\Vert v^{2}-1\Vert U_{3}+\Vert1\Vert U_{1}.
\]
Using this relation we have
\[
\int v^{2}f\,dv=\sum_{ij}X_{i}S_{ij}\left(\Vert v^{2}-1\Vert\langle U_{3},V_{j}\rangle_{v}+\Vert1\Vert\langle U_{1},V_{j}\rangle_{v}\right)=\Vert v^{2}-1\Vert K_{3}+\Vert1\Vert K_{1}.
\]
Thus, for the energy density we have
\[
e=\tfrac{1}{2}\Vert v^{2}-1\Vert K_{3}+\tfrac{1}{2}\Vert1\Vert K_{1}+\tfrac{1}{2}E^{2}.
\]

Our goal is once again to derive the corresponding continuity equation.
We have
\begin{align*}
\partial_{t}e & =\tfrac{1}{2}\Vert v^{2}-1\Vert\partial_{t}K_{3}+\tfrac{1}{2}\Vert1\Vert\partial_{t}K_{1}+E\cdot(\partial_{t}E)\\
 & =\tfrac{1}{2}\Vert v^{2}-1\Vert(U_{3},D[f])_{v}+\tfrac{1}{2}\Vert1\Vert(U_{1},D[f])_{v}+E\cdot(\partial_{t}E)\\
 & =\frac{1}{2}\int v^{2}D[f]\,dv+E\cdot(\partial_{t}E).
\end{align*}
Evaluating the integral and using the relations for the electric field
derived in section \ref{sec:vp-conservation} we get the desired continuity
equation
\[
\partial_{t}e+\nabla_{x}\cdot Q=E\cdot(\partial_{t}E-j)
\]
which implies conservation of energy. 

We have here outlined an approach with $m=3$, $U_{1}\propto1$, $U_{2}\propto v$,
and $U_{3}\propto v^{2}-1$ which conserves mass, momentum, and energy.
However, it is also possible to use this method to only conserve one
or two out of these three invariants. If only energy conservation
is desired, e.g., we can simply choose $U_{1}\propto v^{2}$ and the
argument above is even simpler as we do not have to take the orthogonality
between $U_{1}$ and $U_{3}$ into account.

We can clearly see that for our algorithm the argument with respect
to conservation made for the Vlasov--Poisson equations in section
\ref{sec:vp-conservation} carries over to the low-rank approximation.
This is by design and highlights the ability of our approach to preserve
the physical structure of the Vlasov--Poisson equations.

\section{Discretization\label{subsec:Time-integration}}

In the previous section we have established that the proposed dynamical
low-rank integrator conserves mass, momentum, and energy. This is
done completely within a continuous formulation. Thus, the only approximation
made is due to the fact that the proposed scheme uses a low-rank representation
of the solution. However, to actually implement the conservative dynamical
low-rank integrator on a computer, we have to also introduce a time
and space discretization. Devising a numerical method that discretizes
the equations of motions for the proposed dynamical low-rank integrator,
while maintaining conservation, is the main purpose of this section. 

In section \ref{subsec:conseuler} we will introduce an explicit integrator
for equations (\ref{eq:evol-X})-(\ref{eq:evol-S}) that preserves
mass and momentum up to machine precision. For dynamical low-rank
approximations integrators that are robust to the presence of small
singular values have been developed recently \citep{LO14}. Unfortunately,
it turns out that the projector splitting integrator can not be easily
adapted to the present situation. In section \ref{subsec:Unconventional-integrator}
we propose a robust dynamical low-rank integrator that fits within
the framework considered in this paper.

\subsection{Conservative Euler scheme\label{subsec:conseuler}}

Simply applying a classic Runge--Kutta scheme to the evolution equations
(\ref{eq:evol-X})-(\ref{eq:evol-S}) does destroy the conservative
properties of the method. As an example, let us consider the classic
explicit Euler scheme
\begin{align*}
S_{kl}^{n+1} & =S_{kl}^{n}+\tau\left(X_{k}^{n}V_{l}^{n},D[f^{n}]\right)_{xv},\\
X_{i}^{n+1} & =X_{i}^{n}+\tau\sum_{k}\left(S^{n}\right)_{ik}^{-1}\biggl[\left(V_{k}^{n},D[f^{n}]\right)_{v}-\sum_{l}X_{l}^{n}\left(X_{l}^{n}V_{k}^{n},D[f^{n}]\right)_{xv}\biggr],\\
W_{p}^{n+1} & =W_{p}^{n}+\tau\sum_{q}((S^{n})^{T}S^{n})_{pq}^{-1}\biggl[\frac{1}{f_{0v}}\sum_{i}S_{iq}^{n}(X_{i}^{n},D[f^{n}])_{x}-\sum_{il}S_{iq}^{n}\left(X_{i}^{n}V_{l}^{n},D[f^{n}]\right)_{xv}V_{l}^{n}\biggr],
\end{align*}
where $\tau$ is the time step size and the upper indices denote the
value of the corresponding quantity at the discrete time $t^{n}$.
The reason why this scheme fails to be conservative is that we use
the inverse of $S^{n}$, i.e. $S$ at time $t^{n}$, to compute $X^{n+1}$.
This inconsistency means that there is no well defined $K_{i}^{n}$
and $K_{i}^{n+1}$ and thus the argument in section \ref{sec:Conservation}
can not be applied. 

However, we can write equation (\ref{eq:evol-X}) in the following
conservative form
\begin{equation}
\partial_{t}\biggl(\sum_{i}X_{i}S_{ik}\biggr)=\left(V_{k},D[f]\right)_{v}.\label{eq:evol-K}
\end{equation}
Applying the explicit Euler scheme to that equation yields
\[
\sum_{i}X_{i}^{n+1}S_{ik}^{n+1}=\sum_{i}X_{i}^{n}S_{ik}^{n}+\tau\left(V_{k}^{n},D[f^{n}]\right)_{v}.
\]

No change is made to the other two equations of motions. Putting this
together we obtain
\begin{align}
S_{kl}^{n+1} & =S_{kl}^{n}+\tau\left(X_{k}^{n}V_{l}^{n},D[f^{n}]\right)_{xv},\label{eq:meuler-S}\\
X_{i}^{n+1} & =\sum_{k}(S^{n+1})_{ik}^{-1}\biggl[\sum_{j}X_{j}^{n}S_{jk}^{n}+\tau\left(V_{k}^{n},D[f^{n}]\right)_{v}\biggr],\label{eq:meuler-X}\\
W_{p}^{n+1} & =W_{p}^{n}+\tau\sum_{qi}((S^{n})^{T}S^{n})_{pq}^{-1}S_{iq}^{n}\biggl[\frac{1}{f_{0v}}(X_{i}^{n},D[f^{n}])_{x}-\sum_{l}\left(X_{i}^{n}V_{l}^{n},D[f^{n}]\right)_{xv}V_{l}^{n}\biggr].\label{eq:meuler-W}
\end{align}

In the following we will call this method the conservative Euler scheme.
This is still a fully explicit method as we can first compute $S^{n+1}$
using equation (\ref{eq:meuler-S}), which is then used in the computation
of $X^{n+1}$. 

The scheme satisfies a discrete version of the continuity equation
\begin{align*}
\frac{\rho^{n+1}-\rho^{n}}{\tau} & =\frac{1}{U_{1}}\frac{K_{1}^{n+1}-K_{1}^{n}}{\tau}=\int_{\Omega_{v}}D[f^{n}]\,dv=-\nabla_{x}\cdot j^{n}
\end{align*}
as can be easily seen from equation (\ref{eq:meuler-X}). Integrating
this equation in $x$ we get $M^{n+1}=M^{n}$ and thus conservation
of mass. 

The conservative Euler scheme also yields the following discrete continuity
equation\textbf{
\[
\frac{j^{n+1}-j^{n}}{\tau}=\Vert v\Vert\frac{K_{2}^{n+1}-K_{2}^{n}}{\tau}=\int_{\Omega_{v}}vD[f^{n}]\,dv=-\nabla_{x}\cdot\sigma^{n}-E^{n}\rho^{n}
\]
}which implies conservation of momentum.

For energy we have
\begin{align*}
\frac{e^{n+1}-e^{n}}{\tau} & =\Vert v^{2}-1\Vert\frac{K_{3}^{n+1}-K_{3}^{n}}{2\tau}+\Vert1\Vert\frac{K_{1}^{n+1}-K_{1}^{n}}{2\tau}+\frac{(E^{n+1})^{2}-(E^{n})^{2}}{2\tau}\\
 & =\frac{1}{2}\int_{\Omega_{v}}v^{2}D[f^{n}]\,dv+E^{n}\frac{E^{n+1}-E^{n}}{\tau}+\frac{(E^{n+1}-E^{n})^{2}}{2\tau}\\
 & =\nabla_{x}\cdot Q^{n}+E^{n}\cdot\left(\frac{E^{n+1}-E^{n}}{\tau}-j^{n}\right)+\frac{(E^{n+1}-E^{n})^{2}}{2\tau}.
\end{align*}
The reason why we have rewritten the electric field in this particular
way should become clear shortly. Integrating the above relation in
space yields
\begin{equation}
\mathcal{E}^{n+1}-\mathcal{E}^{n}=\tau\int_{\Omega_{x}}E^{n}\cdot\left(\frac{E^{n+1}-E^{n}}{\tau}-j^{n}\right)\,dx+\frac{1}{2}\int_{\Omega_{x}}(E^{n+1}-E^{n})^{2}dx.\label{eq:energy-diff}
\end{equation}

We again write the electric field using its potential, i.e. $E^{n}=-\nabla\phi^{n}$,
to get
\begin{align*}
\int_{\Omega_{x}}E^{n}\cdot\left(\frac{E^{n+1}-E^{n}}{\tau}-j^{n}\right)\,dx & =\int_{\Omega_{x}}\phi^{n}\left(\frac{\nabla\cdot E^{n+1}-\nabla\cdot E^{n}}{\tau}-\nabla\cdot j^{n}\right)\,dx\\
 & =-\int_{\Omega_{x}}\phi^{n}\left(\frac{\rho^{n+1}-\rho^{n}}{\tau}+\nabla\cdot j^{n}\right)\,dx=0.
\end{align*}
Thus, the first in term on the right-hand side of (\ref{eq:energy-diff})
vanishes. That leaves the second term
\[
\int_{\Omega_{x}}(E^{n+1}-E^{n})^{2}\,dx=\mathcal{O}(\tau^{2}),
\]
which can clearly not be zero (except for trivial solutions that satisfy
$E^{n}=\text{const}$ for all $n$). Thus, the conservative Euler
scheme commits a first order error in the energy. The reason for this
is that in the discrete setting the relation $\partial_{t}(E^{2})=2E\partial_{t}E$
does not hold true. This introduces the second term in the calculations
above and thus destroys the conservation of energy in the discrete
setting.

Fundamentally, the issue is that the Euler scheme is not symmetric
under time reversal. If one allows for implicit methods this deficiency
can be remedied. For example, using $E^{n+1/2}=(E^{n+1}+E^{n})/2$
in the kinetic update results in
\begin{align*}
\frac{e^{n+1}-e^{n}}{\tau} & =\nabla_{x}\cdot Q^{n}-E^{n+1/2}\cdot j^{n}+\frac{(E^{n+1}-E^{n})(E^{n+1}+E^{n})}{2\tau}\\
 & =\nabla_{x}\cdot Q^{n}+E^{n+1/2}\cdot\left(\frac{E^{n+1}-E^{n}}{\tau}-j^{n}\right).
\end{align*}
Integrating in physical space, as above, then shows conservation of
energy. However, since the electric field depends on the particle-density
function this approach results in a fully implicit scheme that has
to be solved up to machine precision, if energy conservation to the
same level of accuracy is desired. We consider the construction of
efficient energy conservative time integrators a subject of future
research.

Let us now turn our attention to the discretization of space. Fortunately,
this is relatively straightforward. As long as the discrete approximation
of the derivatives and the quadrature rule chosen to define the invariants
allows us to perform integration by parts without introducing any
error, the calculations made in this section carry over to the fully
discretized case. This is true for a number of space discretization
strategies. For example, using fast Fourier techniques (FFT) or the
standard second-order centered finite difference scheme to compute
the derivatives in combination with the trapezoidal rule to evaluate
the integrals satisfies this property. The former will be used in
the numerical results in section \ref{sec:numerical-experiments}.
Moreover, for a number of more advanced and higher-order finite difference,
finite volume, and discontinuous Galerkin schemes this property is
satisfied as well (see, e.g., \citep{arakawa1997,einkemmer2014conservative,wiesenberger2019reproducibility}).
\textcolor{black}{Such methods can then also be used in situations
where, e.g.~, homogeneous Dirichlet boundary conditions are required.
For non-homogeneous boundary conditions additional modifications have
to be made to the low-rank integrator (see, e.g., \citep{kusch21}).}

There is one additional concern in an actual implementation that,
by necessity, operates in finite precision arithmetic (i.e.~with
doubles or floats on a computer). To show conservation we employed
the orthogonality condition $\langle W_{p},U_{a}\rangle=0.$ This
is true for the dynamical low-rank integrator by construction (see
section \ref{sec:cons-dlr}). The conservative Euler discretization
also preserves this property. This can be seen explicitly from equation
(\ref{eq:meuler-W}) by taking the inner product with $U_{a}$

\begin{align*}
\langle W_{p}^{n+1},U_{a}\rangle_{v} & =\langle W_{p}^{n},U_{a}\rangle_{v}+\tau\sum_{qi}((S^{n})^{T}S^{n})_{pq}^{-1}S_{iq}^{n}\biggl\langle\frac{1}{f_{0v}}(X_{i}^{n},D[f^{n}])_{x}-\sum_{l}\left(X_{k}^{n}V_{l}^{n},D[f^{n}]\right)_{xv}V_{l}^{n},U_{a}\biggr\rangle_{v}.
\end{align*}
Since $\langle W_{p}^{0},U_{a}\rangle_{v}=0$ is true for the initial
value we can assume that $\langle W_{p}^{n},U_{a}\rangle_{v}=0$ and
thus
\[
\langle W_{p}^{n+1},U_{a}\rangle_{v}=\tau\sum_{qi}((S^{n})^{T}S^{n})_{pq}^{-1}S_{iq}^{n}\left[(X_{i}^{n}U_{a},D[f^{n}])_{xv}-(X_{k}^{n}U_{a},D[f^{n}])_{xv}\right]=0.
\]

However, in the course of the time integration round-off errors can
accumulate and it can thus still happen that the $W_{p}$ fail to
be exactly orthogonal to the $U_{a}$. This can cause a small linear
drift in the error of mass and momentum. To avoid this we perform
an orthonormalization of the $W_{p}$ at the end of each step. This
does not negatively impact the accuracy of the numerical method and
the incurred computational cost is negligible.

\subsection{Unconventional integrator\label{subsec:Unconventional-integrator}}

The integrator described in the previous section is not robust if
the matrix $S$ has small singular values. The reason for this is
that inverting the matrix $S$ is then numerically ill-conditioned.
Small singular values commonly occur if the rank of the solution is
lower than the rank chosen to conduct the numerical simulation. We
refer to \citep{Kieri2016} for a more detailed discussion, but note
that robustness is generally considered a desirable property especially
if the rank is adaptively changed during the simulation.

For the classic dynamical low-rank approximation a projector splitting
integrator has been proposed by Lubich \& Oseledets \citep{LO14}
that remedies this deficiency. The method has also been extended to
a variety of tensor formats \citep{lubich15tio,Lubich2017,El18,ceruti2020time1,ceruti2020time2}.
The main utility of the projector splitting integrator is that by
using a QR decomposition it avoids the inversion of $S$. Unfortunately,
this projector splitting integrator can not be used in the present
situation, because if we solve equation (\ref{eq:evol-W}) to obtain
$\sum_{ip}S_{iq}S_{ip}W_{p}$, while holding the $X_{i}$ constant,
it is not possible to use a QR decomposition to extract the low-rank
factors $S$ and $W$. 

Instead of the projector splitting integrator, we consider an approach
based on the the recently developed unconventional integrator by Ceruti
\& Lubich \citep{ceruti2020unconventional}. The main idea is that
the dynamical low-rank approximation has two parts. On the one hand,
the low-rank factors $X_{i}$ and $V_{j}$ determine the subspaces
$\overline{X}=\text{span}\{X_{i}\}$ and $\overline{V}=\text{span}\{V_{j}\}$
in which an approximation is sought. However, it does not matter how
exactly the $X_{i}$ are chosen as long as they span the same approximation
space $\overline{X}$. On the other hand, the low-rank factor $S$
contains the coefficients that combine the basis functions in an appropriate
way in order to obtain a good approximation. Inverting $S$ is only
required in equations (\ref{eq:evol-X}) and (\ref{eq:evol-W}). That
is, it is only required to obtain the low-rank factors $X_{i}$ and
$V_{j}$.

We thus proceed as follows. First, we discretize equation (\ref{eq:evol-K})
to compute
\[
K_{k}^{n+1}=K_{k}^{n}+\tau\left(V_{k}^{n},D[f^{n}]\right)_{v}\qquad\text{with}\qquad K_{k}=\sum_{i}X_{i}S_{ik}.
\]
and discretize equation (\ref{eq:evol-W}) to compute 
\[
L_{q}^{n+1}=L_{q}^{n}+\frac{\tau}{f_{0v}}\sum_{i}S_{iq}^{n}(X_{i}^{n},D[f^{n}])_{x}-\tau\sum_{il}S_{iq}^{n}\left(X_{k}^{n}V_{l}^{n},D[f^{n}]\right)V_{l}^{n}
\]

with
\[
L_{q}^{n+1}=\sum_{ip}S_{iq}^{n}S_{ip}^{n}W_{p}^{n+1}\qquad\text{and}\qquad L_{q}^{n}=\sum_{ip}S_{iq}^{n}S_{ip}^{n}W_{p}^{n}.
\]

 The approximation spaces $\overline{V}^{n+1}$ and $\overline{X}^{n+1}$
are uniquely defined by $L_{q}^{n+1}$ and $K_{k}^{n+1}$, respectively.
Thus, we perform a QR decomposition 
\[
L_{q}^{n+1}=\sum_{p}W_{p}^{n+1}R_{pq}^{2},\qquad K_{k}^{n+1}=\sum_{i}X_{i}^{n+1}R_{ik}^{1}
\]
to obtain $W_{p}^{n+1}$ and $X_{i}^{n+1}$. The $R$ parts of the
QR decomposition, i.e. $R_{pq}^{2}$ and $R_{ik}^{1}$, are simply
discarded. We then determine $S^{n+1}$ as follows

\begin{align*}
S_{kl}^{n+1} & =\sum_{ij}M_{ki}S_{ij}^{n}N_{jl}^{T}+\tau\left(X_{k}^{n+1}V_{l}^{n+1},D[f(X^{n+1},MS^{n}N^{T},V^{n+1})]\right)_{xv},
\end{align*}
where
\[
M_{ki}=\langle X_{k}^{n+1}X_{i}^{n}\rangle_{x},\qquad N_{jl}^{T}=\langle V_{j}^{n}V_{l}^{n+1}\rangle_{v}.
\]

This procedure is a discretization of equation (\ref{eq:evol-S}),
where we have used the fact that
\[
\sum_{ij}X_{i}^{n}S_{ij}^{n}V_{j}^{n}\approx\sum_{kijl}X_{k}^{n+1}M_{ki}S_{ij}^{n}N_{jl}^{T}V_{l}^{n+1}.
\]
That is, we transform the coefficient matrix $S_{ij}^{n}$ to the
new basis spanned by the previously determined $V_{j}^{n+1}$ and
$X_{i}^{n+1}$ and then solve equation (\ref{eq:evol-S}) with the
basis held fixed over one time step. This gives us a robust dynamical
low-rank integrator as the matrix $S$ needs not be inverted. The
downside of this integrator, however, is that mass and momentum conservation
up to machine precision is lost. \textcolor{black}{The reason for this
is that 
\[
\sum_{i}M_{ki}X_{k}^{n+1},\qquad\qquad\text{and}\qquad\qquad\sum_{j}N_{jl}^{T}V_{l}^{n+1}
\]
}

\textcolor{black}{project the subspaces $\overline{X}^{n}$ and $\overline{V}^{n}$
onto the subspaces $\overline{X}^{n+1}$ and $\overline{V}^{n+1}$,
respectively. Since this projection is not exact, conservation can
be lost in the process. We will also see this in the numerical experiments
conducted in section \ref{sec:numerical-experiments}.}

\section{Numerical experiments\label{sec:numerical-experiments}}

In this section we illustrate the conservative dynamical low-rank
scheme using a number of numerical examples. The numerical algorithm
outlined in sections \ref{sec:cons-dlr} and \ref{subsec:Time-integration}
will be employed. In order to compute the derivatives in space we
employ techniques based on the fast Fourier transform (FFT). 

\subsection{Landau damping\label{subsec:llpert}}

As the first example we consider a variant of Landau damping on the
domain $(x,v)\in[0,4\pi]\times[-6,6]$. The initial value is given
by

\begin{equation}
f(0,x,v)=\left(1+\alpha\cos(\tfrac{1}{2}x)\right)\frac{e^{-v^{2}/2}}{\sqrt{2\pi}}+\sum_{k=1}^{5}\epsilon_{k}\cos(\tfrac{k+1}{2}x)v^{k}e^{-v^{2}/2},\label{eq:llcons}
\end{equation}
where $\alpha=10^{-2}$, $\epsilon_{k}=10^{-4}$ for $k\in\{1,2,3\}$
and $\epsilon_{k}=10^{-5}$ for $k\in\{4,5\}$. This is the classic
Landau damping problem with an added perturbation that ensures that
the initial value is rank $6$.

\begin{figure}[H]
\begin{centering}
\includegraphics[width=14cm]{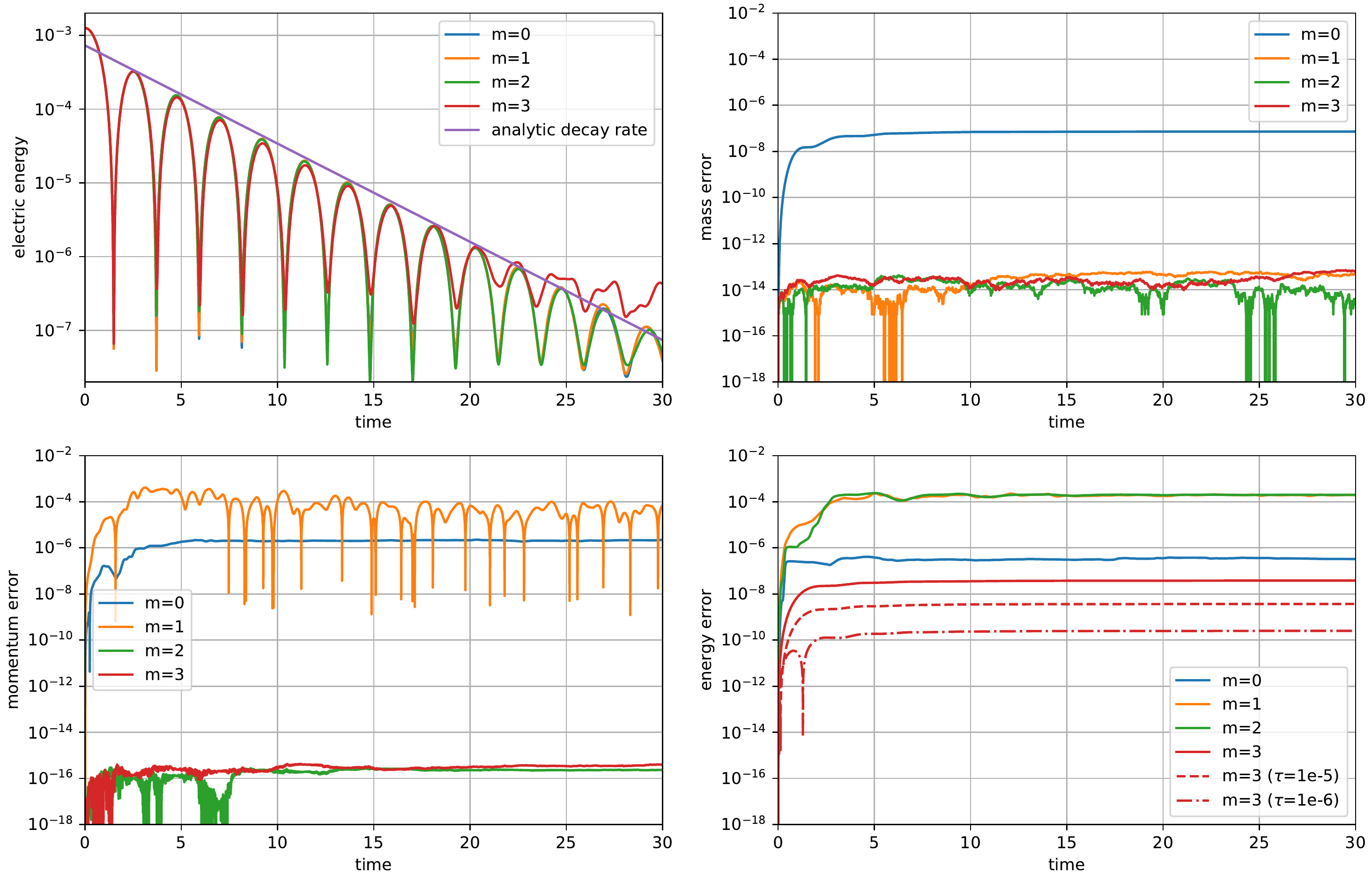}
\par\end{centering}
\caption{Time evolution of the electric energy (top-left), error in mass (top-right),
momentum (bottom-left), and energy (bottom-right) for the linear Landau
damping problem with initial value (\ref{eq:llcons}) is shown. For
mass and energy the relative error and for momentum the absolute error
is reported. The simulation is conducted using the conservative Euler
dynamical low-rank integrator with step size $\tau=10^{-4}$ (except
if otherwise indicated), rank $r=6$, and $128$ grid points in both
the spatial and velocity direction are used. The configurations considered
are $m=0$, $m=1$ ($U_{1}\propto1$), $m=2$ ($U_{1}\propto1$, $U_{2}\propto v$),
and $m=3$ ($U_{1}\propto1$, $U_{2}\propto v$, $U_{3}\propto v^{2}-1$).
\label{fig:llcons}}
\end{figure}

The results of a numerical simulation with the conservative Euler
dynamical low-rank integrator and rank $r=6$ is shown in Figure \ref{fig:llcons}.
For this problem the linearized decay rate can be determined analytically
and matches the observed numerical results well. To study the conservation
properties we run simulations using $m=0$ (no conservation), $m=1$
($U_{1}\propto1$ and thus mass conservation), $m=2$ ($U_{1}\propto1$,
$U_{2}\propto v$ and thus mass and momentum conservation), and $m=3$
($U_{1}\propto1$, $U_{2}\propto v$, $U_{3}\propto v^{2}-1$ and
thus mass, and momentum conservation as well as energy conservation
in absence of any time integration error). The numerical results in
Figure \ref{fig:llcons} show clearly that mass and momentum are conserved
up to machine precision. Thus, the results match perfectly with the
theory in sections \ref{sec:Conservation} and \ref{subsec:Time-integration}.
Energy is interesting as mandating only mass or momentum conservation
can increase the error in energy. This is not entirely unexpected
as similar behavior can be observed even for the case of Hamiltonian
ordinary differential equations, see e.g.~ \citep[Example 4.3 in Chap. IV.4][]{hairer2006geometric}.
However, one should note that for $m=3$ the error in energy is significantly
reduced compared to all the other configurations. Thus, using the
energy conservative dynamical low-rank integrator is clearly beneficial,
even if the time integration error is taken into account. Reducing
the time step size further also improves energy conservation, indicating
that the energy error in the $m=3$ configuration is limited by the
time integration error, as expected. We also note that for the $m=3$
configuration the fidelity of the simulation is somewhat worse. The
reason for this that the dynamical low-rank integrator is only able
to choose $3$ basis functions (all others are already determined
by imposing conservation) and has thus less freedom to decrease the
error in the particle density function (see also the discussion in
section \ref{subsec:numerical-ll}). 

\subsection{Maxwellian with position-dependent velocity}

As our second example, we consider a Maxwellian particle density

\begin{align}
f(0,x,v) & =\frac{\rho(x)}{\sqrt{2\pi}}\exp\left(-(v-u(x))^{2}/2\right)\label{eq:f0before}
\end{align}
with a position-dependent velocity $u(x)=\alpha\cos(x)$ on the domain
$(x,v)\in[0,4\pi]\times[-6,6]$. The density is given by $\rho=1+\epsilon\cos(\tfrac{1}{2}x)$.
In the simulation the parameters are chosen as $\alpha=0.2$ and $\epsilon=10^{-2}$.
Note that the initial value given in (\ref{eq:f0before}) is not low-rank.
However, for small to moderate $u$, as is the case here, we can use
the following expansion (see, e.g., \citep{chen1998lattice,E18})

\begin{align*}
\exp\left(-(v-u(x))^{2}/2\right)=\exp(-v^{2}/2) & \left(1+uv+\tfrac{1}{2}u^{2}(v^{2}-1)+\tfrac{1}{6}u^{3}(v^{3}-3v)+\tfrac{1}{24}u^{4}(3-6v^{2}+v^{4})\right.\\
 & \left.\;+\tfrac{1}{120}u^{5}(15v-10v^{3}+v^{5})\right)+\mathcal{O}(u^{5}).
\end{align*}

Thus, for the initial value of the dynamical low-rank integrator we
use the following rank $6$ approximation

\begin{equation} \label{eq:llnonconstu} 
\begin{split}
f(0,x,v)=\frac{\rho}{\sqrt{2\pi}}\exp(-v^{2}/2)&\left(1+uv+\tfrac{1}{2}u^{2}(v^{2}-1)+\tfrac{1}{6}u^{3}(v^{3}-3v)+\tfrac{1}{24}u^{4}(3-6v^{2}+v^{4})\right.
\\
\left.\;+\tfrac{1}{120}u^{5}(15v-10v^{3}+v^{5})\right)
\end{split}
\end{equation}

\begin{figure}[H]
\begin{centering}
\includegraphics[width=14cm]{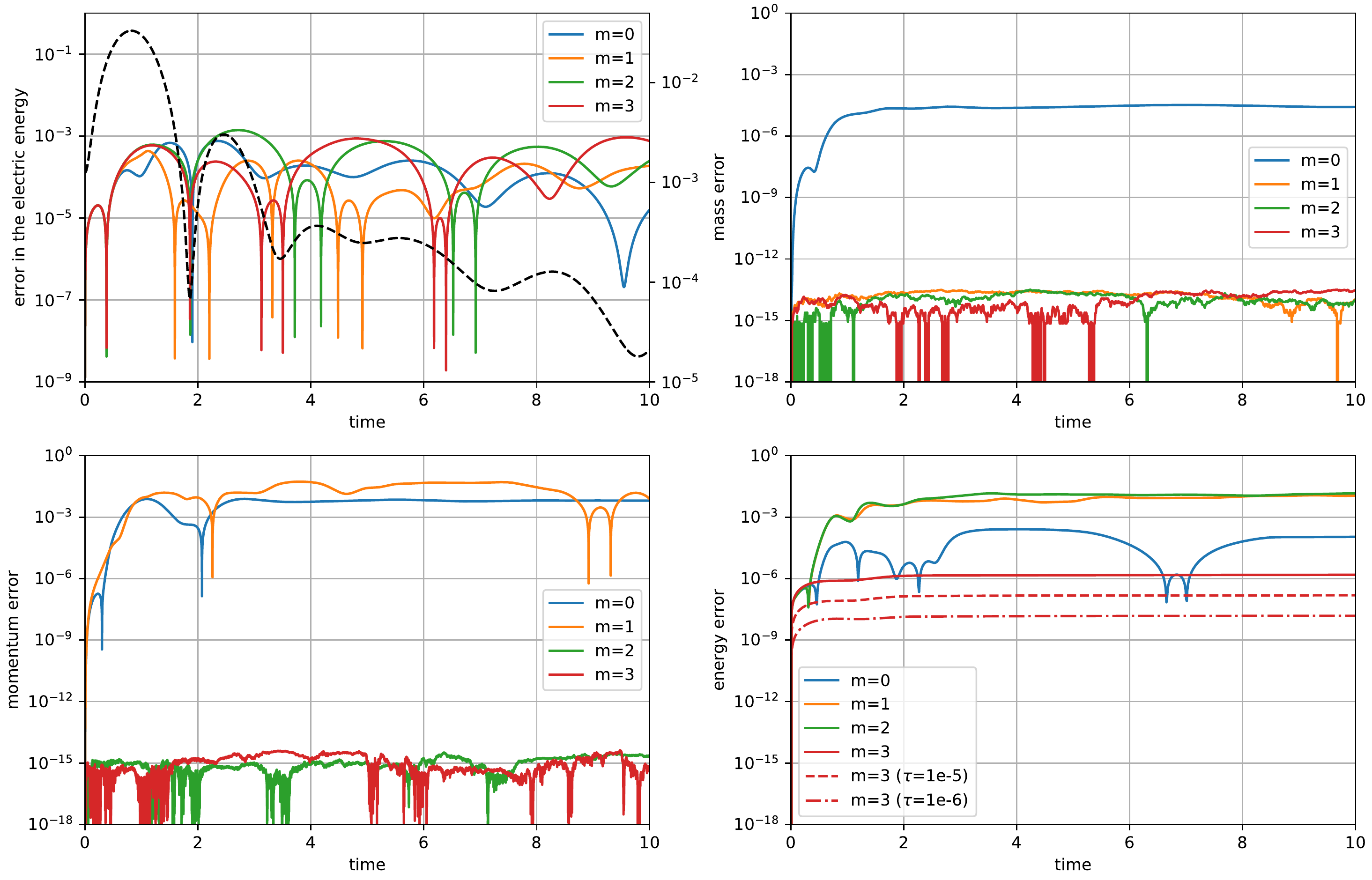}
\par\end{centering}
\caption{Time evolution of the error in electric energy (top-left), error in
mass (top-right), error in momentum (bottom-left), and error in energy
(bottom-right) for the Maxwellian with position-dependent velocity
(\ref{eq:f0before}) is shown. To determine the error in the electric
energy a reference solution is computed using a full grid semi-Lagrangian
scheme that uses $1024$ degrees of freedom in both the $x$ and $v$
direction. For mass and energy the relative error and for momentum
the absolute error is reported. The simulation is conducted using
the conservative Euler dynamical low-rank integrator with step size
$\tau=10^{-4}$ (except if otherwise indicated), rank $r=6$, and
$128$ grid points in both the spatial and velocity direction are
used. The configurations considered are $m=0$, $m=1$ ($U_{1}\propto1$),
$m=2$ ($U_{1}\propto1$, $U_{2}\propto v$), and $m=3$ ($U_{1}\propto1$,
$U_{2}\propto v$, $U_{3}\propto v^{2}-1$). \textcolor{black}{The
time evolution of the electric energy is shown with a dashed black
line (axis on the right).} \label{fig:llnonconstu}}
\end{figure}

The results of a numerical simulation with the conservative Euler
dynamical low-rank integrator and rank $r=6$ is shown in Figure \ref{fig:llnonconstu}.
We again observe excellent agreement between theory and the numerical
results. In particular, mass and momentum are conserved up to machine
precision (and are between $9$ and $14$ orders of magnitude smaller
than for the $m=0$ configuration) and the error in energy is reduced
and shown to be dominated by the time integration error.

\subsection{Landau damping using the unconventional dynamical low-rank integrator\label{subsec:numerical-ll}}

\begin{figure}[H]
\begin{centering}
\includegraphics[width=14cm]{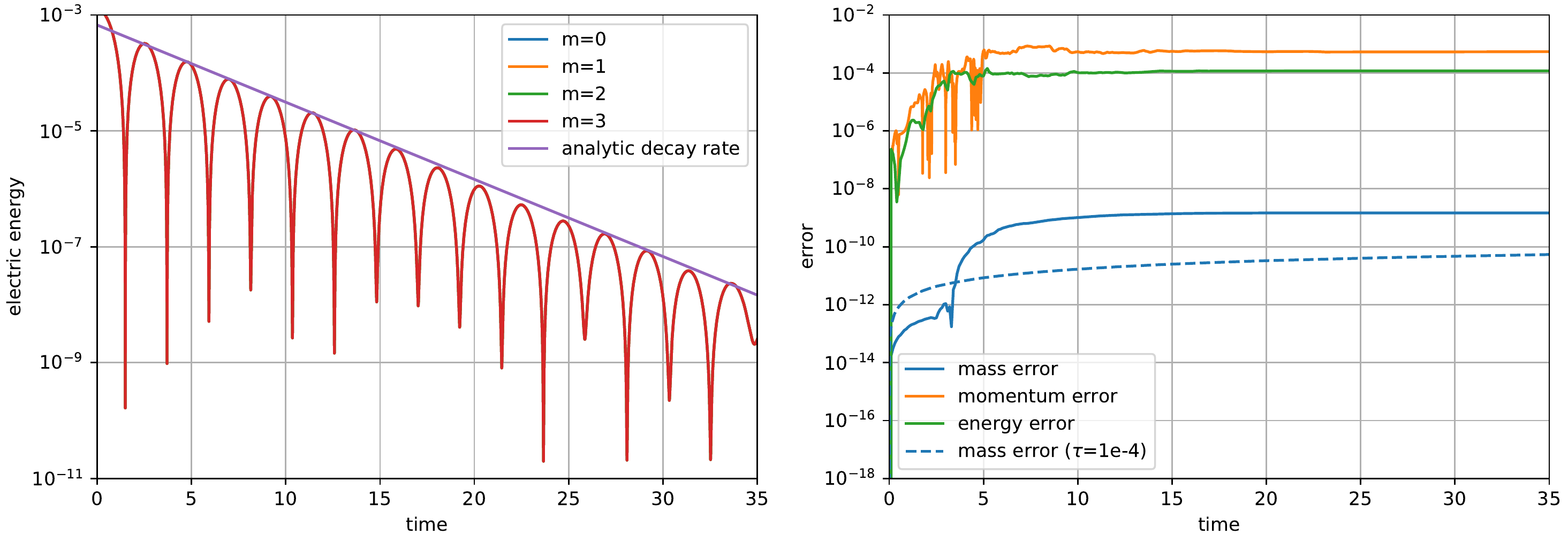}
\par\end{centering}
\caption{Time evolution of the electric energy (left) and error in mass, momentum,
and energy (right) for the linear Landau damping problem with initial
value (\ref{eq:ll}) is shown. For mass and energy the relative error
and for momentum the absolute error is reported. The simulation is
conducted using the unconventional dynamical low-rank integrator with
step size $\tau=10^{-3}$ , rank $r=10$, and $128$ grid points in
both the spatial and velocity direction are used. The configurations
considered are $m=0$, $m=1$ ($U_{1}\propto1$), $m=2$ ($U_{1}\propto1$,
$U_{2}\propto v$), and $m=3$ ($U_{1}\propto1$, $U_{2}\propto v$,
$U_{3}\propto v^{2}-1$). The results for all configurations nearly
overlap in the plot. \label{fig:ll}}
\end{figure}

The conservative Euler dynamical low-rank integrator can not be used
if the rank with which the simulation is run is larger than the rank
of the solution (see the discussion in section \ref{subsec:Unconventional-integrator}).
This is often the case for a number of commonly considered plasma
instabilities, where the initial value has rank $1$. However, for
the unconventional dynamical low-rank integrator (see section \ref{subsec:Unconventional-integrator})
this is not an issue. To demonstrate this we will consider the following
linear Landau damping problem
\begin{equation}
f(0,x,v)=\left(1+\alpha\cos(\tfrac{1}{2}x)\right)\frac{e^{-v^{2}/2}}{\sqrt{2\pi}},\label{eq:ll}
\end{equation}

which is rank $1$. For the simulation we use $\alpha=10^{-2}$ on
the domain $(x,v)\in[0,4\pi]\times[-6,6]$. The results of the numerical
simulation with rank $r=10$ is shown in Figure \ref{fig:ll}. We
observe that the analytic decay rate is reproduced accurately by all
configurations. In particular, no reduction in fidelity is observed
for the $m=3$ configuration (as is the case in section \ref{subsec:llpert}).
Thus, having a certain number of basis functions that the algorithm
can choose freely is clearly beneficial. As explained in section \ref{subsec:Unconventional-integrator},
the unconventional integrator is not conservative up to machine precision.
The error in mass, momentum, and energy is dominated by the time integration
error and thus reducing the time step size reduces the error in the
conserved quantities, as is illustrated for mass in Figure \ref{fig:ll}. 

\subsection{Two-stream instability using the unconventional dynamical low-rank
integrator}

\begin{figure}[H]
\begin{centering}
\includegraphics[width=14cm]{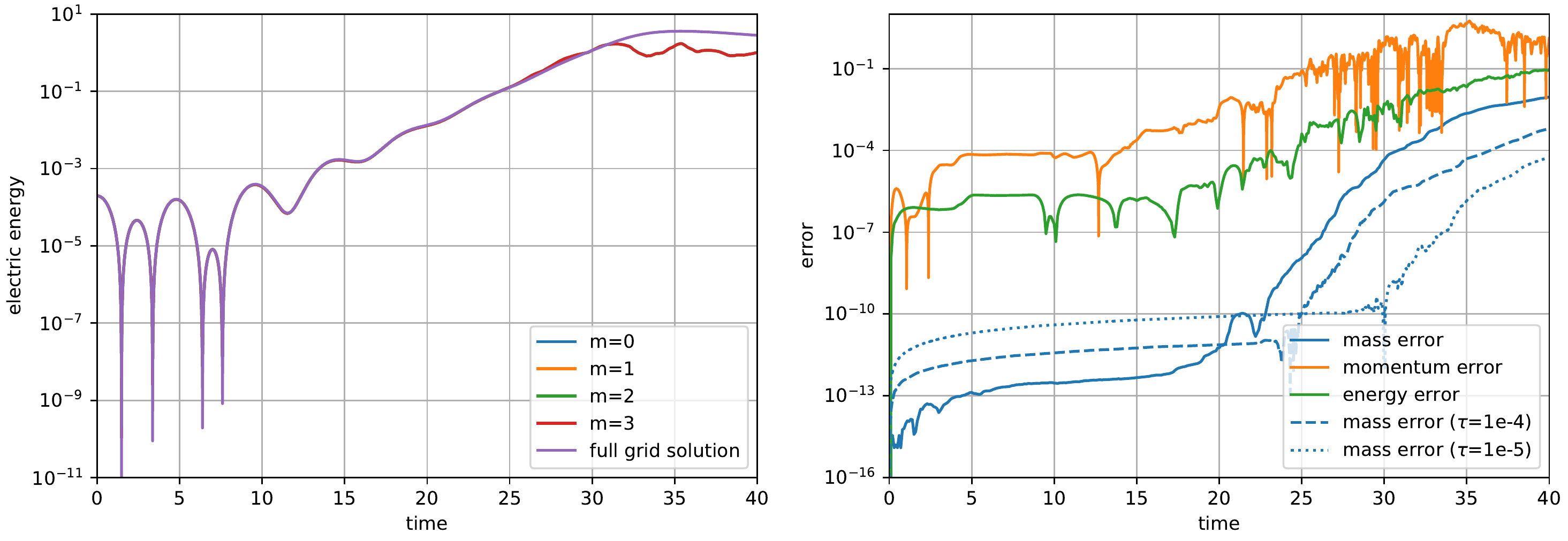}
\par\end{centering}
\caption{Time evolution of the electric energy (left) and error in mass, momentum,
and energy (right) for the two-stream instability (\ref{eq:two-stream})
is shown. For mass and energy the relative error and for momentum
the absolute error is reported. The simulation is conducted using
the unconventional dynamical low-rank integrator with step size $\tau=10^{-3}$
, rank $r=10$, and $128$ grid points in both the spatial and velocity
direction are used. The configurations considered are $m=0$, $m=1$
($U_{1}\propto1$), $m=2$ ($U_{1}\propto1$, $U_{2}\propto v$),
and $m=3$ ($U_{1}\propto1$, $U_{2}\propto v$, $U_{3}\propto v^{2}-1$).
The results for all configurations nearly overlap in the plot. The
full grid solution is computed using a semi-Lagrangian scheme that
uses $1024$ degrees of freedom in both the $x$ and $v$ direction.
\textcolor{black}{Note that initially the error in mass is dominated
by round-off and thus smaller time step sizes lead to a larger error.}
\label{fig:two-stream}}
\end{figure}

Finally, we consider the two-stream instability
\begin{equation}
f(0,x,v)=\tfrac{1}{2\sqrt{2\pi}}\left(e^{-(v^{2}-\bar{v})/2}+e^{-(v^{2}+\bar{v})/2}\right)\left(1+\alpha\cos(kx)\right)\label{eq:two-stream}
\end{equation}
with $\alpha=10^{-3}$, $k=0.2$, and $\bar{v}=2.4$ on the domain
$(x,v)\in[0,10\pi]\times[-7,7]$. This initial value has rank $1$.
However, the rank of the solution increases significantly as the system
is evolved in time and nonlinear effects become stronger. Nevertheless,
it is known that the dynamical low-rank integrator resolves the corresponding
dynamics well, at least, up to saturation \citep{El18,EL18_cons}.
From the numerical results in Figure \ref{fig:two-stream} we see
that this is also the case for the conservative dynamical low-rank
integrator. As in the previous section the error in mass, momentum,
and energy is dominated by the time integration error in all cases.

\section*{Acknowledgements}

Work by I. Joseph was performed under the auspices of the U.S. DOE
by LLNL under Contract DE-AC52-07NA27344.

\bibliographystyle{plain}
\bibliography{references}

\end{document}